\documentstyle[11pt]{article}

\input amssym.tex 

\font\tinybbfont=msbm6
\font\scriptsizebbfont=msbm7 scaled \magstep 1
 1
\font\footnotesizebbfont=msbm9 scaled \magstep 0
\font\smallbbfont=msbm7 scaled \magstep 2
\font\bbfont=msbm9 scaled \magstep1  
 1
\font\Largebbfont=msbm10 scaled \magstep 2
 3
 4

\def\tinyBbb#1{\hbox{\tinybbfont #1}}
\def\scriptsizeBbb#1{\hbox{\scriptsizebbfont #1}}

\def\footnotesizeBbb#1{\hbox{\footnotesizebbfont #1}}
\def\smallBbb#1{\hbox{\smallbbfont #1}}
\def\Bbb#1{\hbox{\bbfont #1}}

\def\LargeBbb#1{\hbox{\Largebbfont #1}}

\newcommand{\CP}{{\Bbb C}{\rm P}}

\newcommand{\GL}{\mbox{\it GL}\,}

\newcommand{\degree}{\mbox{{\it deg}\,}}
\newcommand{\ev}{\mbox{{\it ev}\,}}

\newcommand{\pt}{\mbox{\it pt}}
\newcommand{\rank}{\mbox{\it rank}\,}

\begin{document}

\enlargethispage{23cm}

\begin{titlepage}

$ $

\vspace{-1cm} 

\noindent\hspace{-1cm}
\parbox{6cm}{March 2000}\
   \hspace{7.5cm}\
   \parbox{5cm}{ {\tt math.AG/0003071} }

\vspace{2cm}     

\centerline{\large\bf
A Reconstruction of Euler Data}

\vspace{2cm}

%
\hspace{-6em}
\begin{minipage}{18cm}
 \begin{center}
  \parbox[t]{5cm}{
   \centerline{\large Bong H.\ Lian$^1$ }
   \vspace{1.1em}
   \centerline{\it Department of Mathematics}
   \centerline{\it Brandeis University}
   \centerline{\it Waltham, MA 02154} } \
  \hspace{1em} \
  \parbox[t]{5cm}{
   \centerline{\large Chien-Hao Liu$^2$}
   \vspace{1.1em}
   \centerline{\it Department of Mathematics}
   \centerline{\it Harvard University}
   \centerline{\it Cambridge, MA 02138}  } \
  \hspace{1em}
  \parbox[t]{5cm}{
   \centerline{\large Shing-Tung Yau$^3$}
   \vspace{1.1em}
   \centerline{\it Department of Mathematics}
   \centerline{\it Harvard University}
   \centerline{\it Cambridge, MA 02138}  }
 \end{center}
\end{minipage}

\vspace{2cm}

\begin{quotation}
\centerline{\bf Abstract}
\vspace{0.3cm}

\baselineskip 12pt  
{\small
We apply the mirror principle of [L-L-Y] to reconstruct the Euler data
$Q=\{Q_d\}_{d\in{\tinyBbb N}\cup\{0\}}$ associated to a vector bundle
$V$ on ${\smallBbb C}{\rm P}^n$ and a multiplicative class $b$.
This gives a direct way to compute the intersection number $K_d$
without referring to any other Euler data linked to $Q$. Here $K_d$
is the integral of the cohomology class $b(V_d)$ of the induced
bundle $V_d$ on a stable map moduli space.
A package '{\tt \verb+EulerData_MP.m+}' in Maple V that carries
out the actual computation is provided.
For $b$ the Chern polynomial, the computation of $K_1$ for the bundle
$V=T_{\ast}{\smallBbb C}{\rm P}^2$, and $K_d$, $d=1,2,3$, for
the bundles ${\cal O}_{{\tinyBbb C}{\rm P}^4}(l)$ with $6\le l\le 10$
done using the code are also included.
} 
\end{quotation}

\bigskip

\baselineskip 12pt
{\footnotesize
 } 

\bigskip

\noindent {\small
MSC number 1991: 14Q15, 68-04, 55N91, 55R91, 81T30.
} 

\bigskip

\baselineskip 11pt
{\footnotesize
\noindent{\bf Acknowledgements.}
We would like to thank
 David Eisenbud, Joe Harris, Yi Hu, Sheldon Katz, Albrecht Klemm,
 Kefeng Liu, Jason Starr, Richard Thomas, Cumrun Vafa, and Eric Zaslow
for valuable discussions at various stages of the work,
 Shinobu Hosono,
who took the extreme pain to go through a preliminary Mathematica
version of the code and gave lots of suggestions for improvement, and to
 Ling-Miao Chou,
who went through and checked the Maple code of its near final version.
Without their advices and involvements, our code may not be born
yet at the moment. C.H.L.\ would like to thank, in addition,
 Hung-Wen Chang and Daniel Freed
for the valuable conversation and discussion on programming before/at
the beginning of the work, and
 Orlando Alvarez, Philip Candelas, Jacques Distler, Martin Halpern,
 Rafael Nepomechie, and Xenia de la Ossa,
who influence and shape his view on the physics side of the subject. 
This work is supported by
 DOE grant DE-FG02-88ER25065 and
 NSF grants DMS-9619884 and DMS-9803347.
} 

\noindent
\underline{\hspace{20em}}

$^1${\footnotesize E-mail: lian@brandeis.edu}

$^2${\footnotesize E-mail: chienliu@math.harvard.edu}

$^3${\footnotesize E-mail: yau@math.harvard.edu}

\end{titlepage}

\newpage
$ $

\vspace{-4em}  

\centerline{\sc  Reconstruction of Euler Data}

\vspace{1.5em}

\baselineskip 14pt  

\begin{flushleft}
{\Large\bf 0. Introduction and outline.}
\end{flushleft}

\begin{flushleft}
{\bf Introduction.}
\end{flushleft}
Ever since the ground-breaking work of 
[C-dlO-G-P], mirror symmetry and its meaning and consequenes have
been investigated by several groups of people both from the
mathematical and from the physical point of view.
(See introduction of [L-L-Y] and references therein for background and
for a comparison of different approaches. See also [MS].)

In this note, we apply the theory developed
in the series of papers [L-L-Y] I, II, and III,
to compute new intersection numbers on stable map moduli spaces.
The theory goes beyond justifying mirror symmmetry that relates
the usually difficult A-model computations to the much more tractable
B-model computations for Calabi-Yau manifolds. Indeed, the theory
is a way to directly do the "A-model"
computations for manifolds which are not necessarily Calabi-Yau. It is the goal
of these notes to explain this and to provide a computer code that
carries out the actual computations for bundles over
$\CP^n$.

To make this article more self-contained, 
we recall in Sec.\ 1 the definitions of the basic objects
involved and the Atiyah-Bott localization formula that is used substantially
in the theory. In Sec.\ 2 and Sec.\ 3, we focus on the case of
critial bundles over $\CP^n$ and consider their Euler classes.
In Sec.\ 2, we give a quick summary of facts and formula from
[L-L-Y] I-III that are directly related to the actual computation
of Euler data $\{Q_d\}_d$ and the intersection numbers $K_d$.
In Sec.\ 3, we explain how the theory of [L-L-Y] gives rise to
a system of linear equations that can be solved inductively.
The solution of the system gives the Euler data $\{Q_d\}_d$, from
which the intersection numbers $K_d$ can be computed.
After these, we then discuss in Sec.\ 4 the modifications needed
to take into account also non-critical bundles. There the Chern
polynomial is considered.
In Sec.\ 5, we single out six examples whose first few $K_d$
are computed this way via a Maple code.
In Sec.\ 6, the Maple code {\tt \verb+EulerData_MP.m+}
with instructions is given. Eighteen cases have been tested and
computed. The last record of the run for each of these cases is
given in SEC.\ 3 of the code for references. The code 
provided can be easily modified to compute
other cases of interest.

This article is served as a supplement to and a computational
account of [L-L-Y]. As a result, our notations and terminologies
follow [L-L-Y] very closely. Readers are referred to ibidem for more
theoretical details.

\bigskip

\begin{flushleft}
{\bf Outline.}
\end{flushleft}
{\small
\baselineskip 11pt  

\begin{quote}
 1. Essential mathematical backgrounds for physicists.

 2. Summary of related constructs in ``Mirror Principle". 

 3. Computation of $Q_d$ inductively.

 4. Modifications for non-critical bundles over
    ${\smallBbb C}{\rm P}^n$.

 5. Examples.

 6. A package in Maple V for the computation of $Q_d$ and $K_d$.
\end{quote}
} 

\baselineskip 14pt  

\newpage

\section{Essential mathematical backgrounds for physicists.}

We collect in this section the most essential backgrounds
for understanding these notes. Along the way, we also set up the
notations for the notes.

\bigskip

\noindent $\bullet$
{\bf Stable maps and their moduli.} [Ko]

\bigskip

\noindent
{\bf Definition 1.1 [stable map].} {Let $X$ be a smooth projective
variety.  An {\it $n$-pointed stable map into $X$} consists
 of a connected marked curve $(C,p_1,\cdots,p_n)$ and a morphism
 $f:C\rightarrow X$ satisfying the following properties:
 \begin{quote}
  \hspace{-1.9em}(1)\hspace{1ex}
  The only singularities of $C$ are ordinary double points.
   
  \hspace{-1.9em}(2)\hspace{1ex}
  $p_1,\,\cdots,\,p_n$ are distinct ordered smooth points of $C$.

  \hspace{-1.9em}(3)\hspace{1ex}
  If $C_i$ is a component of $C$ that is isomorphic $\CP^1$ and is
mapped to a point under $f$, then $C_i$ contains at least three
  special (i.e.\ nodal or marked) points.
 
  \hspace{-1.9em}(4)\hspace{1ex}
  If $C$ has (arithmertic) genus $1$ and $n=0$, then $f$ is not
  constant.
 \end{quote}

\bigskip

\noindent
{\it Remark 1.2.} Given Conditions (1) and (2) in the above
 definition, Conditions (3) and (4) are equivalent to the assertion
 that the data $(f,C,p_1,\,\cdots,\,p_n)$ has only finitely many
 automorphisms.
} 

\bigskip

Given a class $\beta\in H_2(X,{\Bbb Z})$, the moduli space of all
stable maps $(f,C,p_1,\cdots,\,p_n)$ such that $[f(C)]=\beta$ into
$X$ will be denoted by $\overline{\cal M}_{g,n}(X,\beta)$.

\bigskip

\noindent $\bullet$
{\bf Equivariant cohomology.} (See [Au].)
Given a group $G$ acting on a space $X$.  Let $BG$ be the classifying
space and $\mbox{\it EG}\rightarrow \mbox{\it BG}$ be the
universal principle $G$-bundle associated to $G$.
The equivariant cohomology $H^{\ast}_G(X)$ of $X$ associated to the
$G$-action is defined to be
$$
 H^{\ast}_G(X)\;=\; H^{\ast}(X_G)\,,
$$
where $X_G = \mbox{\it EG}\times_G X$ is the total space of the
associated $X$-bundle over $BG$. Note that
$H^{\ast}_G(\pt)=H^{\ast}(BG)$, and that 
$H^{\ast}_G(X)$ is naturally a  $H^{\ast}_G(\pt)$-module.

The constant map $X\rightarrow \pt$ induces an equivariant projection
$\pi_X:X_G\rightarrow BG$. The induced pushforward map
${\pi_X}_!$ from $H^{\ast}_G(X)$ to $H^{\ast}_G(\pt)$ is given by
integration along the fiber of $\pi_X$. This is also called the 
{\it equivariant integral}. In notation,
$$
 {\pi_X}_!\;=\;\int_{X_G}\;:\;
             H^{\ast}_G(X)\;\longrightarrow\; H^{\ast}_G(\pt)\,.
$$

\bigskip

\noindent
{\it Remark 1.3.}
In this article, the coefficient for $H^{\ast}(X_G)$ can be
${\Bbb Q}$, ${\Bbb R}$, or ${\Bbb C}$. Usually we use ${\Bbb Q}$ or
${\Bbb C}$ in the discusssion.

\bigskip

\noindent
{\bf Example 1.4.}
Let ${\Bbb T}^r=\prod_rS^1$ be an $r$-torus. Then
$B{\Bbb T^r}=\prod_r\CP^{\infty}$ and
$H^{\ast}(B{\Bbb T}^r)=H^{\ast}_{{\scriptsizeBbb T}^r}(\pt)
                        ={\Bbb C}[\lambda_1,\,\cdots,\,\lambda_r]$,
the polynomial ring generated by $\lambda_1,\,\cdots,\,\lambda_r$,
where $\lambda_i$ is the first Chern class of the hyperplane line
bundle ${\cal O}(1)$ over the $i^{th}$ $\CP^{\infty}$ in the
product.

Let ${\Bbb T}^r\rightarrow \GL(N+1,{\Bbb C})$ be a representation
of an $r$-torus on ${\Bbb C}^{N+1}$ with weight
$(\beta_0,\,\cdots,\,\beta_N)$. Note that each $\beta_i$ is a
linear combination of the $\lambda_j$'s. This induces a
${\Bbb T}^r$-action on $\CP^N$. With respect to this action,
$$
 H^{\ast}_{{\scriptsizeBbb T}^r}(\CP^N)\;
 =\; \left. H^{\ast}_{{\scriptsizeBbb T}^r}(\pt)[\zeta]
      \mbox{\raisebox{-2ex}{\rule{0ex}{1ex}}} \right/
    \mbox{\raisebox{-.4ex}{$\prod_{i=0}^N(\zeta-\beta_i)$}}\,,
$$
where $\zeta$ is the equivariant hyperplane class from a fixed
lifting of the hyperplane class of $\CP^N$. The equivariant integral
$\int_{({\scriptsizeBbb C}{\rm P}^n)_{{\scriptsizeBbb T}^r}}:\,
      H^{\ast}_{{\scriptsizeBbb T}^r}(\CP^N)\,
             \rightarrow\, H^{\ast}_{{\scriptsizeBbb T}^r}(\pt)$
picks out the coefficient of $\zeta^N$ in elements of
$H^{\ast}_{{\scriptsizeBbb T}^r}(\CP^N)$.

\noindent\hspace{14cm} $\Box$

\bigskip

\noindent $\bullet$
{\bf The Atiyah-Bott localization formula.} [A-B]
Let $T$ be an $r$-torus that acts on a manifold $X$ with the
set of fixed points a union of smooth connected submanifolds $Z_j$.
Then the normal bundle $N_j$ of $Z_j$ in $X$ is a $T$-equivariant
vector bundle with its equivariant Euler class
$e_T(N_j)\in H^{\ast}_T(Z_j)$.

There are three fundamental maps between $H^{\ast}_T(Z_j)$ and
$H^{\ast}_T(X)\,$:
\begin{quote}
  \hspace{-1.9em}(1)\hspace{1ex}
  {\it the restriction homomorphism}$\,$:
  $$
   i^{\ast}_j\; :\; H^{\ast}_T(X)\;\longrightarrow\;H^{\ast}_T(Z_j)
  $$
  induced by the equivariant inclusion $i_j:Z_j\hookrightarrow X$;

  \hspace{-1.9em}(2)\hspace{1ex}
  {\it the Gysin map}$\,$:
  $$
   {i_j}_!\;:\; H^{\ast}_T(Z_j)\; \longrightarrow\; H^{\ast}_T(X)\,;
  $$
\end{quote}
Note that ([Au]), for any $\alpha_j\in H^{\ast}_T(Z_j)$, one has
$$
 i^{\ast}_j\,\circ\,{i_j}_!(\alpha)\;=\; \alpha_j\,\cup\,e_T(N_j)\,.
$$

Let ${\cal R}$ be the localization of $H^{\ast}_T(\mbox{\it BT})$.
With the notation following Example 1.4,
${\cal R}={\Bbb C}(\lambda_1,\,\cdots,\,\lambda_r)$. Then
$e_T(N_j)$ is an invertible element in the
localization $H^{\ast}_T(Z_j)\otimes{\cal R}$. We can now state
the Atiyah-Bott localization formula [A-B]:

\bigskip

\noindent
{\bf Fact 1.5 [Atiyah-Bott localization formula].} {\it
 The following map is an isomorphism
 $$
  \begin{array}{ccccl}
    H^{\ast}_T(X)\otimes{\cal R} 
      & \stackrel{\sim}{\longrightarrow}
      & \bigoplus_j\,H^{\ast}_T(Z_j)\otimes{\cal R} & \\[.6ex]
    \alpha & \longmapsto
      & \left(\, \frac{i^{\ast}_j(\alpha)}{e_T(N_j)}\,\right)_j &.
  \end{array}
 $$
 Its inverse is given by
 $$
  (\alpha_j)_j\;\longmapsto\; \sum_j\, {i_j}_!(\alpha_j)\,.
 $$
 Combining these two, one has
 $$
  \alpha\; =\; \sum_j\, {i_j}_!\,
            \left(\,\frac{i^{\ast}_j(\alpha)}{e_T(N_j)}\,\right)\,,
 $$
 for any $\alpha\in H^{\ast}_T(X)\otimes{\cal R}$.
} 

\bigskip

\noindent
{\bf Example 1.6.} ({\it Continuing Example 1.4}).
The fixed point set on $\CP^N$ comes from the $N+1$ coordinate lines
of ${\Bbb C}^{N+1}$. Denote this set by $\{\,p_0,\,\cdots,\,p_N\,\}$.
Let $\alpha\in H^{\ast}_T(\CP^N)\otimes{\cal R}$, then $\alpha$ can be
written as a polynomial $f(\zeta)$ with coefficients in ${\cal R}$.
In terms of this, $i^{\ast}_j(\alpha)=f(\beta_j)$;
$e_T(N_j)=\prod_{k\ne j}(\beta_j-\beta_k)$, where $k$ runs in
$\{\,0,\,\cdots,\,N\,\}$; and the Gysin map ${i_j}_!$ is given by
the cup product with $\prod_{k\ne j}(\zeta-\beta_k)$.
Thus, from the localization formula, one has 
$$
 f(\zeta)\; =\; \sum_{j=0}^N\,f(\beta_i)\,
                   \frac{\prod_{k\ne j}(\zeta-\beta_k)}{
                           \prod_{k\ne j}(\beta_j-\beta_k)}\,.
$$
\hspace{14cm} $\Box$

\bigskip

\noindent $\bullet$
{\bf Concavex bundles over $\CP^n$.} ([L-L-Y, I].)
Let $T={\Bbb T}^{n+1}$ acts on ${\Bbb C}^{n+1}$ with weights
$\lambda_0,\,\cdots,\,\lambda_n$. It induces an action on $\CP^n$.

\bigskip

\noindent
{\bf Definition 1.7 [concavex bundle].}{
 Let $V$ be a $T$-equivariant vector bundle over $\CP^n$. We call
 $V$ {\it convex} (resp.\ {\it concave}) if the $T$-equivariant
 Euler class $e_T(V)$ is invertible and $H^1(C,f^{\ast}V)=0$
 (resp.\ $H^0(C,f^{\ast}V)=0$) for any $0$-pointed genus $0$ stable
 map $f:C\rightarrow \CP^n$. We call $V$ {\it concavex} if it is a
 direct sum of a convex and a concave bundle. We denote this
 decomposition by $V=V^+\oplus V^-$.
}

\bigskip

\noindent
{\bf Definition 1.8 [splitting type].} {
Let $V$ be a $T$-equivariant concavex bundle over $\CP^n$.
Let 
$l_a$, $k_b$ be positive integers such that for every
$T$-invariant line $C\cong\CP^1$ in $\CP^n$ we have
a $T$-equivariant isomorphism
$$
 V|_C\;\cong\;\oplus_a\,{\cal O}(l_a)\,\oplus\,\oplus_b\,{\cal O}(-k_b).
$$
Then we call
$(\,l_1,\,\cdots;\,k_1,\,\cdots\,)$ the {\it splitting type} of $V$.
}

\bigskip

\section{Summary of related constructs in ``$\,$Mirror Principle$\,$".}

This section follows [L-L-Y, I]. Readers should consult [L-L-Y, I] 
(see also [L-L-Y, II]) for more  details.

\bigskip

\noindent $\bullet$
{\bf Set-up of the problem.}
Let $T={\Bbb T}^{n+1}$, $V$ be a concavex $T$-equivariant bundle
over $\CP^n$ of splitting type $(\,l_1,\,\cdots;\,k_1,\,\cdots\,)$,
and $d$ be a class in $H_2(\CP^n,{\Bbb Z})={\Bbb Z}$. Let
$\overline{\cal M}_{0,0}(\CP^n,d)$ be the moduli space of stable
maps into $\CP^n$, of degree $d$, genus $0$, without marked points;
and similarly for $\overline{\cal M}_{0,1}(\CP^n,d)$ and
$M_d=\overline{\cal M}_{0,0}(\CP^1\times\CP^n,(1,d))$.
Let $N_d=\CP^{(n+1)d+n}$ be the linear $\sigma$-model for $\CP^n$.
(We will say more about $N_d$ in the next item.) Then one has the
following $S^1\times T$-equivariant diagram:
$$
 \begin{array}{cccccccccl}
  & & \hspace{-3em}V_d\,=\,\pi^{\ast}U_d &
    & \hspace{-.4em}U_d & & \hspace{-.4em}\rho^{\ast}U_d & & V &\\[.6ex]
  & & \hspace{-4em}\downarrow & & \hspace{-1em}\downarrow &
                                & \downarrow  & & \downarrow & \\[.6ex]
  N_d  & \stackrel{\varphi}{\longleftarrow}   & \hspace{-4em}M_d
       & \hspace{-6em}\stackrel{\pi}{\longrightarrow}
       & \hspace{-2em}\overline{\cal M}_{0,0}(\CP^n,d)
       & \stackrel{\rho}{\longleftarrow}
       & \overline{\cal M}_{0,1}(\CP^n,d)
       & \stackrel{\mbox{\footnotesize\it ev}}{\longrightarrow}
       & \CP^n         \\[.6ex]
  & &  \hspace{-4em}\|     &&&&&&  \\[.6ex]
  & &  \hspace{-4em}\overline{\cal M}_{0,0}(\CP^1\times\CP^n,(1,d))
               &&&&&&&,  
 \end{array}
$$
where $\rho$ forgets and $\ev$ evaluates at the marked point of a
$1$-pointed stable map; $U_d\,=\,\rho_!\,\ev^{\ast}\,V$,
the pushforward via $\rho$ of the pullback of $V$ via $\ev$;
$\pi$, the {\it contracting morphism}, is induced by the projection
of a stable map in $\CP^1\times\CP^n$ to the $\CP^n$ component and
contracting the unstable components;
and $\varphi$, {\it the collapsing morphism}, will be explained in
the next two items.

Let $c_{\mbox{\scriptsize\it top}}$ be the top Chern class
(i.e.\ the Euler class) of $U_d$, then the
{\it intersection number of degree $d$} is defined to be
$$
 K_d\; =\;
  \int_{\overline{\cal M}_{0,0}({\scriptsizeBbb C}{\rm P}^n,d)}\,
  c_{\mbox{\scriptsize\it top}}(U_d)\,.
$$
One of the goals in the mirror symmetry literatures is to compute
$K_d$'s and to relate them to enumerative problems on $\CP^n$, or
some varieties therein. 
An important insight from [L-L-Y] is that one can
reduce this problem to an easy problem on projective spaces
$\{N_d\}_{d=0}^{\infty}$, 
called the linear $\sigma$-model in [L-L-Y].
In fact the intersection numbers $K_d$ can be recovered from
cohomology classes
$Q=\{Q_d\}_{d=0}^{\infty}$, called Euler data, defined
on those projective spaces. In turn the Euler data
can be computed essentially by an elementary  algorithm,
and, sometimes, via an explicit formula.
A nonlinear recursion involving graph sums was
use to compute $K_4$ in the case of $O(5)$ on $\CP^4$
in [Ko].

If $b$ is any multiplicative cohomology classes,
then we can apply our algorithm to compute the integrals
$$
 K_d\; =\;
  \int_{\overline{\cal M}_{0,0}({\scriptsizeBbb C}{\rm P}^n,d)}\,
  b(U_d)\,.
$$
More generally, suppose $a,\,b,\,\ldots$ are any multiplicative
cohomology classes. Then, for any given vector bundle $V$, we have
$$
 a(V)=a_0(V)+a_1(V)+\cdots+a_r(V)\,,
$$
where $a_i(V)$ is the degree $i$ component of $a(V)$.
We can homogenize $a$ by writing
$$
 a_x(V)=\sum_i x^{r-i} a_i(V)\,,
$$
where $x$ is a formal variable, and the class $a_x$ remains
multiplicative. Likewise, we have
$$
 b_y(V)=\sum_i y^{r-i} b_i(V)\,,
$$
etc. Multiplying them, we get $a_x(V) b_y(V)\cdots$, which is a new
multiplicative class. Note that this is a polynomial in the variables
$x,\,y,\,\ldots$, with coefficients of the form $a_i(V)b_j(V)\cdots$.
Moreover, each such product correspond to a unique monomial in
$x,\,y,\,\ldots$. Now our algorithm computes
$K_d\in{\Bbb C}[x,\,y,\,\ldots]$ for the multiplicative class
$a_x(U_d) b_y(U_d)\cdots\,$ and, hence, computes all coefficients of
the form
$$
  \int_{\overline{\cal M}_{0,0}({\scriptsizeBbb C}{\rm P}^n,d)}\,
   a_i(U_d)b_j(U_d)\,, \ldots\,.
$$
Note that this number is zero unless the integrand has the right
total degree.

\bigskip

\noindent
{\it Notation.} (Cf.\ Example 1.4.) We shall adapt the
following notations for the rest of the article:
 $G=S^1\times T$,
 $\lambda=(\lambda_0,\,\cdots,\,\lambda_n)$,
 $\alpha=c_1({\cal O}(1))\in H^{\ast}_{S^1}(BS^1)$,
 ${\cal R}={\Bbb Q}(\lambda)[\alpha]$,
 ${\cal R}^{-1}={\Bbb Q}(\lambda,\alpha)$,
 ${\cal R}H^{\ast}_G(\,\cdot\,)=H^{\ast}_G(\,\cdot\,)
       \otimes_{{\scriptsizeBbb Q}[\lambda,\alpha]}{\cal R}$,
 ${\cal R}^{-1}H^{\ast}_G(\,\cdot\,)=H^{\ast}_G(\,\cdot\,)
       \otimes_{{\scriptsizeBbb Q}[\lambda,\alpha]}{\cal R}^{-1}$.

\bigskip

\noindent $\bullet$
{\bf The linear $\sigma$-model $\{N_d\}_{d=0}^{\infty}$ for
     ${\CP^n}$.}
Let
$$
 N_d\; =\; {\Bbb P}\,(\,H^0(\CP^1,
      {\cal O}_{{\scriptsizeBbb C}{\rm P}^1}(d))^{n+1})\;
       \cong\; \CP^{(n+1)d+n}
$$
be the space of $(n+1)$-tuple of homogeneous polynomials of degree
$d$ on $\CP^1$ up to an overall constant multiple in ${\Bbb C}$.
An element in $N_d$ can be written as
$$
 [\,\sum_r z_{0r}w_0^rw_1^{d-r}\,:
            \,\cdots\,:\,\sum_r z_{nr}w_0^r w_1^{d-r}\,]\,,
$$
where $[w_0:w_1]$ is the homogeneous coordinates for $\CP^1$ and
$z_{ir}\in{\Bbb C}$. The sequence $\{N_d\}_{d=0}^{\infty}$ is called
{\it the linear $\sigma$-model} for $\CP^n$.

Let $G=S^1\times{\Bbb T}^{n+1}$; then $G$ acts on
$\CP^1\times\CP^n$ by
$$
 (t,t_0,\,\cdots,\,t_n)\,
     \cdot\,([w_0:w_1]\,,\,[x_0:\,\cdots\,:x_n])\;
 =\; ([tw_0:w_1]\,,\,[t_0x_0:\,\cdots\,:t_nx_n])\,.
$$
This induces a $G$-action on $N_d$ with fixed points
$$
 p_{i,r}\;=\; [\,0\,:\,\cdots\,:\,0\,:\,
                w_0^rw_1^{d-r}\,:\,0\,:\,\cdots\,:\,0\,]\,,
$$
where the non-zero term appears at the $i^{th}$ position,
$i=0,\,\cdots,\,n$, and $r=0,\,\cdots,d$. Note that the weight for
the $G$-action at $T_{p_{i,r}}N_d$ is $\lambda_i+r\alpha$.

There are two $G$-equivariant maps between the $N_d$'s, defined
as follows:
\begin{quote}
 $I:N_{d-1}\rightarrow N_d\,$, \hspace{2ex} 
 $[f_0:\cdots:f_n]\mapsto[w_1f_0:\cdots:w_1f_n]\,$; \hspace{2ex} and

 $\overline{\mbox{\rule{0em}{1.2ex}$\hspace{1.5ex}$}}:
   N_d\rightarrow N_d\,$, \hspace{2ex}
 $[f_0(w_0,w_1):\cdots:f_n(w_0,w_1)]
         \mapsto[f_0(w_1,w_0):\cdots:f_n(w_1,w_0)]$.
\end{quote}
From $I$, one obtaines a chain of inclusions
$$
 N_0=\CP^n\;\stackrel{I}{\longrightarrow}\; N_1\;
  \stackrel{I}{\longrightarrow}\;\cdots\;
  \stackrel{I}{\longrightarrow}\;N_d\,,
$$
whose composition gives a canonical incluion
$I_d:N_0=\CP^n\rightarrow N_d$. Let $\kappa$ be the equivariant
hyperplane class in $H^{\ast}_G(N_d)$. Then the induced map of
$\overline{\mbox{\rule{0em}{1.2ex}$\hspace{1.5ex}$}}$ on
${\cal R}^{-1}H^{\ast}_G(N_d)$ is generated by
$\overline{\kappa}=\kappa-d\alpha$, $\overline{\alpha}=-\alpha$,
and $\overline{\lambda_i}=\lambda_i$.

\bigskip

\noindent $\bullet$
{\bf The $G$-equivariant morphism $\varphi:M_d\rightarrow N_d$.}
First note that $M_d$ and $N_d$ are two different compactifications
of the space ${\cal M}_{0,0}(\CP^n,d)$ of degree $d$ maps from
$\CP^1$ to $\CP^n$. Precisely, an element in
${\cal M}_{0,0}(\CP^n,d)$ can be written as
$$
 [f_0:\,\cdots\,:f_n]\;
    =\;[\,\sum_r z_{0r}w_0^rw_1^{d-r}\,:
            \,\cdots\,:\,\sum_r z_{nr}w_0^r w_1^{d-r}\,]
$$
with $f_0,\,\cdots,\,f_n$ relatively prime. Its embedding in $N_d$
is tautological, while its embedding in $M_d$ is given by
$$
 [f_0:\,\cdots\,:f_n]\;\longmapsto\;
  ([w_1:w_0]\,,\,[f_0:\,\cdots\,:f_n])\,.
$$
Via these embeddings, the identity map on ${\cal M}_{0,0}(\CP^n,d)$
extends to a $G$-equivariant morphism $\varphi:M_d\rightarrow N_d$.
Explicitly, $\varphi$ can be described as follows:
\begin{quote}
 Let $(f,C)\in M_d$ and $\pi_1$, $\pi_2$ be the projections of
 $\CP^1\times\CP^n$ onto its first and second factor respectively.
 Then one can decompose $C$ into $C_0\cup C_1\cup\,\cdots\,\cup C_s$
 with $C_0\cap C_j=x_j$ for $j>0$ such that
 $\pi_1\circ f:C_0\,\stackrel{\sim}{\rightarrow}\,\CP^1$
 and any other $C_j$ is pinched to some 
 $\pi_1\circ f(x_j)=[a_j,b_j]\in\CP^1$ under $\pi_1\circ f$.
 Let $d_i$ be the degree of $\pi_2\circ f:C_j\rightarrow \CP^n$ and
 $[\sigma_0:\,\cdots\,:\sigma_n]$ represent the degree $d_0$ map
 $\pi_2\circ f:C_0\rightarrow \CP^n$. Then
 $$
  \varphi\;:\, (f,C)\;\longmapsto\;
                            [g\,\sigma_0:\,\cdots\,:g\,\sigma_n]\,,
   \hspace{1em}\mbox{where}\hspace{1em}
   g\;=\; \prod_{j=1}^s\,(a_jw_0-b_jw_1)^{d_j}\,.
 $$
\end{quote}

\bigskip

\noindent $\bullet$
{\bf Euler data.}

\bigskip

\noindent
{\bf Definition 2.1 [Euler data].} {
 Given an invertible class $\Omega\in H^{\ast}_T(\CP^n)^{-1}$,
 the localization of $H^{\ast}_T(\CP^n)$, an
 {\it $\Omega$-Euler data} is a sequence $Q=\{Q_d\}_{d=0}^{\infty}$
 of classes $Q_d\in{\cal R}H^{\ast}_G(N_d)$ that satisfy
 \begin{quote}
  \hspace{-1.9em}(1)\hspace{1ex}
   $Q_0=\Omega$.

  \hspace{-1.9em}(2)\hspace{1ex} The {\it gluing identity}$\,$:
   $$
    i^{\ast}_{p_i}(\Omega)\,i^{\ast}_{p_{i,r}}(Q_d)\;
     =\; \overline{i^{\ast}_{p_{i,0}}(Q_r)}\:
            i^{\ast}_{p_{i,0}}(Q_{d-r})\,,
   $$
   for all $d$ and $i=0,\,\cdots,\,n$, $r=0\,\cdots\,d$.
 \end{quote}
} 

\bigskip

An immediate consequence is the following lemma:

\bigskip

\noindent
{\bf Fact 2.2 [reciprocity].} (Lemma 2.4 in [L-L-Y, I].) {\it
 If $Q$ is an Euler data, then, for $i,j=0,\,\cdots,\,n$,
 $r=0,\,\cdots,\,d$, $d=0,\,1,\,2,\,\cdots$, one has
 \begin{quote}
  \hspace{-1.9em}(1)\hspace{1ex}
  $Q_d(\lambda_i+d\alpha)=\overline{Q_d(\lambda_i)}$.

  \hspace{-1.9em}(2)\hspace{1ex}
  $Q_d(\lambda_i)|_{\alpha=(\lambda_i-\lambda_j)/d}\;
         =\; Q_d(\lambda_j)|_{\alpha=(\lambda_j-\lambda_i)/d}$
  for $d>0$.

  \hspace{-1.9em}(3)\hspace{1ex}
  $\Omega(\lambda_i)\,Q_d(\lambda_j)\,
                           =\,Q_r(\lambda_j)Q_{d-r}(\lambda_i)$
  at $\alpha=(\lambda_j-\lambda_i)/\mbox{\raisebox{-.4ex}{$r$}}$
  for $r>0$.
 \end{quote}
} 

\bigskip

Recall the various bundles and maps from Item Set-up. 

\bigskip

\noindent
{\bf Fact 2.3 [Euler data].} (Theorem 2.8 in [L-L-Y, I].) {\it
 Let $V=V^+\oplus V^-$ be a concavex bundle over ${\CP^n}$,
 $\chi^V_d$ be the equivariant Euler class of $V_d$,
 $Q_0=\Omega^V=e_T(V^+)/\mbox{\raisebox{-.4ex}{$e_T(V^-)$}}$,
 $Q_d=\varphi_!(\chi^V_d)$ for $d>0$. Then 
 $Q=\{Q_d\}$ is an $\Omega^V$-Euler data.
} 

\bigskip

We call a concavex bundle $V\rightarrow\CP^n$ {\it critical} if the
induced bundle $U_d\rightarrow\overline{\cal M}_{0,0}(\CP^n,d)$ has
rank equal to
$\mbox{\it dim}\,(\overline{\cal M}_{0,0}(\CP^n,d))=(n+1)d+n-3$.

\bigskip

\noindent
{\bf Fact 2.4.} (Theorem 3.2 (ii) in [L-L-Y, I].)
{\it
 Let $V$ be a critical concavex bundle over $\CP^n$. Then in the
 non-equivariant limit $\lambda\rightarrow 0$,
 $$
  \int_{{\scriptsizeBbb C}{\rm P}^n}
   e^{-Ht/{\alpha}}\,
    \frac{ \lim_{\lambda\rightarrow 0} I^{\ast}_d(Q_d)}{
               \prod_{m=1}^d(H-m\alpha)^{n+1}}\;
  =\; \alpha^{-3}\,(2-dt)\,K_d\,.
 $$
} 

\bigskip

\noindent
Thus, once $Q_d$ is determined, the intersection number $K_d$ is
also determined.

\bigskip

\noindent
{\it Remark 2.5.}
Since $H^{n+1}=0$, one can rewrite the above formula as
{\small
$$
 \int_{{\scriptsizeBbb C}{\rm P}^n}
  \left[\,\sum_{k=0}^n\frac{(-Ht/\alpha)^k}{k!}\,\right]\,
  \left[\, \lim_{\lambda\rightarrow 0}I^{\ast}_d(Q_d)\,\right] \,
  \frac{(-1)^{(n+1)d}}{(d!)^{n+1}\alpha^{(n+1)d}}\,
  \left[\, \sum_{k=0}^n
               \left( \frac{H}{m\alpha} \right)^k\, \right]^{n+1}\;
 =\; \alpha^{-3}\,(2-dt)\,K_d\,.
$$
{\normalsize Note} } 
also that in [C-K], there is another formula, implicitly in [L-L-Y, I],
that relates $Q_d$ and $K_d$:
$$
 \int_{{\scriptsizeBbb C}{\rm P}^n}
  H\, e^{-Ht/{\alpha}}\,
   \frac{ \lim_{\lambda\rightarrow 0} I^{\ast}_d(Q_d)}{
              \prod_{m=1}^d(H-m\alpha)^{n+1}}\;
 =\; \alpha^{-2}\,d\,K_d\,.
$$

\bigskip

\noindent $\bullet$
{\bf Determination of an Euler data.}
By the localization formula, $Q_d$ is determined by its restriction
$i^{\ast}_{p_{i,r}}(Q_d)$ at the fixed points $p_{i,r}$, for
$i=0,\,\cdots,\,n$, $r=0,\,\cdots,\,d$.
Explicitly,
$$
 Q_d\;=\; \sum_{(i,r)}\,
         \frac{i^{\ast}_{p_{i,r}}(Q_d)\,
           \prod_{(j,s)\ne(i,r)}(\kappa-\lambda_j-s\alpha)}{
             \prod_{(j,s)\ne(i,r)}
         \left(\,\rule{0ex}{2ex}
                  \lambda_i-\lambda_j+(r-s)\alpha\,\right)}\,.
$$
Since the Euler data condition says that 
$$
 i^{\ast}_{p_{i,r}}(Q_d)\;
  =\; \frac{\overline{i^{\ast}_{p_{i,0}}(Q_r)}\;
       i^{\ast}_{p_{i,0}}(Q_{d-r})}{i^{\ast}_{p_i}(\Omega)}\,,
$$
it turns out that, to determine $Q_d$, one only needs to know its
restrictions $i^{\ast}_{p_{i,0}}(Q_d)$ at $p_{i,0}$ for
$i=0,\,\cdots,\,n$.

We can now state theorems from [L-L-Y, I] that enables one to
determine $i^{\ast}_{p_{i,0}}(Q_d)$.

\bigskip

\noindent
{\bf Fact 2.6 [degree bound and determination of Euler data].} 
 (Theorem 2.10, Theorem 2.11, and Theorem 3.2 (i) in [L-L-Y, I].)
 {\it 
 Let $V$ be a concavex bundle over $\CP^n$ of splitting type
 $(l_1,l_2,\cdots; k_1,k_2,\cdots)$, $\chi$ be its Euler
 characteristic class, and
 $Q=\{Q_d\}_{d\in{\scriptsizeBbb N}\cup\{0\}}$ be the $\chi$-Euler
 data for $V$, as in Fact 2.3.  Then the restrictions
 $I^{\ast}_d(Q_d)\in H^{\ast}_G(\CP^n)$ has
 $$
  \degree_{\alpha}I^{\ast}_d(Q_d)\;\le \;(n+1)d-2\,.
 $$
 Furthermore, $Q$ is completely determined by the value of the
 restrictions
 $i^{\ast}_{p_{i,0}}(Q_d)$, $i=0,\cdots,n$, $d=0,1,2,\cdots$,
 at $\alpha=(\lambda_i-\lambda_j)/\mbox{\raisebox{-.4ex}{$d$}}$,
 $i\ne j$. These values are given explicitly by
 $$
  \left. i^{\ast}_{p_{i,0}}(Q_d)
   \right|_{\alpha=\frac{\lambda_i-\lambda_j}{d}}\;
  =\; \prod_a\prod_{m=0}^{l_ad}\,
               (l_a\lambda_i-m\,\frac{\lambda_i-\lambda_j}{d})\:
      \prod_b\prod_{m=1}^{k_bd-1}\,
               (-k_b\lambda_i+m\,\frac{\lambda_i-\lambda_j}{d})\,.
 $$
} 

\bigskip

\noindent
{\it Remark 2.7 [total degree bound].}
For $V$ critical, since the rank of
$U_d\rightarrow \overline{\cal M}_{0,0}(\CP^n, d)$
is $(n+1)d+n-3$, the total degree of $Q_d$, as a polynomial of
$\kappa$, $\alpha$, and $\lambda$, is bounded by $(n+1)d+n-3$. 

\bigskip

These facts and remarks allow one to compute $Q_d$ as a polynomial
of $\kappa$ and $\alpha$ with coefficients in
${\Bbb C}(\lambda_0,\,\cdots,\,\lambda_n)$. We now turn to this
detail.

\bigskip

\section{Computation of $Q_d$ inductively.}

Following previous notations, let $V$ be a critical concavex bundle
over $\CP^n$ of splitting type
$(l_1,\, l_2,\,\cdots;\, k_1,\,k_2,\,\cdots\,)$.
Then $Q_1$ can be computed, using the 
Atiyah-Bott formula. The higher $Q_d$
can be computed by the recursive relation from the gluing identities,
the special values of $Q_d$ at the fixed points, and the
$\alpha$ degree bound of $I^{\ast}_d(Q_d)$.

\bigskip

\begin{flushleft}
{\bf The computation of $Q_1$.}
\end{flushleft}
For $d=1$, $\degree_{\alpha}Q_1(\lambda_i,\alpha)\le n-1$
and, for a fixed $i$, the $n$-many values
$Q_1(\lambda_i,\lambda_i-\lambda_j)$, $j\ne i$
are known from Fact 2.6:
$$
 Q_1(\lambda_i, \lambda_i-\lambda_j)\;
   =\; \prod_a\prod_{m=0}^{l_a}\,
                (l_a\lambda_i-m\,(\lambda_i-\lambda_j))\:
       \prod_b\prod_{m=1}^{k_b-1}\,
                (-k_b\lambda_i+m\,(\lambda_i-\lambda_j))\,,
        \hspace{1em}\mbox{for}\hspace{1.5ex}j\ne i\,.
$$
Thus, using the Lagrange interpolation formula, one obtains
$$
 i^{\ast}_{p_{i,0}}(Q_1)\;=\; Q_1(\lambda_i,\alpha)\;
 =\; \sum_{ \scriptsize
           \begin{array}{c}j=0,\,\cdots,\,n \\ j\ne i\end{array} }
      Q_1(\lambda_i,\lambda_i-\lambda_j)\,
      \frac{ \prod_{k\ne i,j}\,(\,\alpha-\lambda_i+\lambda_k\,) }{
         \prod_{k\ne i,j}\,(\,\lambda_k-\lambda_j\,)}\,.
$$
By the Reciprocity Lemma, 
$i^{\ast}_{p_{i,1}}(Q_1)\,=\,\overline{Q_1(\lambda_i,\alpha)}\,
  =\,Q_1(\lambda_i,-\alpha)$.
In this way, the restriction of $Q_1$ at the set of fixed points
of the $S^1\times {\Bbb T}^{n+1}$-action on $N_1$ are all acquired.
Using the localization formula and playing around with the indices,
one obtains an exact expression
{\small
$$
 Q_1\; =\;
  \sum_{i=0}^n\, \left[\,
     \left(\, \rule{0ex}{2ex}
                 f_i(\alpha)\,(\kappa-\lambda_i-\alpha)\,
                 +\, f_i(-\alpha)\,(\kappa-\lambda_i)\, \right)\,
             \prod_{j\ne i}(\kappa-\lambda_j)\,
              \prod_{j\ne i}(\kappa-\lambda_j-\alpha)\,\right]\,,
$$
{\normalsize where}
} 
{\small
\begin{eqnarray*}
 \lefteqn{
    f_i(\alpha)\;=\;\frac{i^{\ast}_{p_{i,0}}(Q_1)}{
         \prod_{(j,s)\ne(i,0)}(\lambda_i-\lambda_j-s\alpha)}}\\[1ex]
 & &  =\; \sum_{j\ne i}\,
   \frac{ \hspace{1ex}
    \prod_a \prod_{m=0}^{l_a}
          \left(\,\rule{0ex}{2ex}
                  l_a\lambda_i-m(\lambda_i-\lambda_j)\,\right)\,
      \prod_b\prod_{m=1}^{k_b-1}
       \left(\rule{0ex}{2ex}
       -k_b\lambda_i+m(\lambda_i-\lambda_j)\,\right) \hspace{1ex} }{
              \alpha\, (\alpha-\lambda_i+\lambda_j)\,
                \prod_{k\ne i}(\lambda_i-\lambda_k)\,
                 \prod_{k\ne i,j} (\lambda_j-\lambda_k) }    \,.
\end{eqnarray*}
} 

\bigskip

\begin{flushleft}
{\bf The computation of $Q_d$ for $d>1$.}
\end{flushleft}
Let $N=(n+1)d+n-3$. Then one may write $Q_d$ as a polynomial in
$\kappa,\alpha$ with coefficients in ${\Bbb C}[\lambda]$:
$$
 Q_d\; =\; \sum_{\mu=0}^N\,\sum_{\nu=0}^{N-\mu}\,
              w_{\mu\nu}\,\alpha^{\mu}\kappa^{\nu}
   \hspace{1em}\mbox{with}\hspace{1em}
   w_{\mu\nu}\in {\Bbb C}[\lambda] \,.
$$
Since we have an explicit formula for $Q_1$, we may assume that
$Q_0,\, Q_1\,\cdots,\,Q_{d-1}$ are all determined.

\bigskip

\noindent $\bullet$
{\it Systems from the gluing identity}$\,$:
The gluing identity that an Euler data must satisfy is completely
encoded in the following two systems of linear equations in
$w_{\mu\nu}$:

\bigskip

\begin{quote}
 \hspace{-1.9em}(1)\hspace{1ex} {\it 
 $Q_d(\lambda_i+r\alpha,\,\alpha)$
 for $i=0,\,\cdots,\,n$ and $r=1,\,\cdots,\, d-1$}$\,$:
 \medskip

 The gluing identity says
 {\small
 $$
  Q_d(\lambda_i+r\alpha,\,\alpha)\;
   =\; \sum_{\mu=0}^N\,\sum_{\nu=0}^{N-\mu}\,
         w_{\mu\nu}\,\alpha^{\mu}(\lambda_i+r\alpha)^{\nu}\;  
   =\; i_{p_{i,r}}(Q_d)\;
   =\;  \left.
       \overline{Q_r(\lambda_i,\alpha)}\,Q_{d-r}(\lambda_i,\alpha)
         \mbox{\raisebox{-1ex}{\rule{0ex}{1ex}}} \right/
     \mbox{\raisebox{-.4ex}{$\Omega(\lambda_i)$}}   
 $$
 {\normalsize for $i=0,\,\cdots,\,n$ and $r=1,\,\cdots,\,d-1$.}}
 Denote
 $\left. \overline{Q_r(\lambda_i,\alpha)}\,Q_{d-r}(\lambda_i,\alpha)
                 \mbox{\raisebox{-1ex}{\rule{0ex}{1ex}}} \right/
                     \mbox{\raisebox{-.4ex}{$\Omega(\lambda_i)$}}$
 by $b_1(i,r)$; then $b_1(i,r)$ is known by induction, Furthermore,
 $$
  \degree_{\alpha}(b_1(i,r))\;
   \le\;  (n+1)r-2+(n+1)(d-r)-2\;=\;(n+1)d-4\; <\; N\,.
 $$
 Thus, one may write $b_1(i,r)=\sum_{s=0}^N\,b_1(i,r,s)\alpha^s$,
 where $b_1(i,r,s)=0$ for $s\ge (n+1)d-3$. After expanding the powers
 and exchanging and relabelling the indices to the above equation,
 one obtains the following linear system in $w_{\mu\nu}\,$:

 \raisebox{-.6em}{$\left\{ \rule{0em}{3em}\right.$} \hspace{1em}
 \parbox{12cm}{
  $$
   \sum_{\mu=0}^s\,\sum_{\nu=s-\mu}^{N-\mu}\,
       w_{\mu\nu}\,\cdot\,
        \left({\small \begin{array}{c}
                        \nu \\ s-\mu
                      \end{array}}\right)\,
          r^{s-\mu}\:\lambda_i^{\mu+\nu-s}\;
       =\; b_1(i,r,s)  \hspace{12em}
  $$
  \hspace{4em} for $i=0,\,\cdots,\,n$, $r=1,\,\cdots,\,d-1$, and
  $s=0,\,\cdots,\,N$.
 } 

 \bigskip

 \hspace{-1.9em}(2)\hspace{1ex} {\it
  $Q_d(\lambda_i+d\alpha,\alpha)$ for $i=0,\,\cdots,\,n$}$\,$: 
  \medskip

  The gluing identity says
  {\small
  $$
   Q_d(\lambda_i+d\alpha,\alpha)\;
    =\; \sum_{\mu=0}^N\,\sum_{\nu=0}^{N-\mu}\,
          w_{\mu\nu}\,\alpha^{\mu}(\lambda_i+d\alpha)^{\nu}\;
    =\;\overline{Q_d(\lambda_i,\alpha)}\;
    =\; Q_d(\lambda_i,-\alpha)\;
    =\; \sum_{\mu=0}^N\,\sum_{\nu=0}^{N-\mu}\,
             w_{\mu\nu}\,(-\alpha)^{\mu}\lambda_i^{\nu}
  $$
  {\normalsize for $i=1=0,\,\cdots,\,n$.}}
  This gives rise to the second linear system in $w_{\mu\nu}\,$:

  \raisebox{-.6em}{$\left\{ \rule{0em}{3em}\right.$}\hspace{1em}
  \parbox{12cm}{
   $$
    \sum_{\mu=0}^{s-1}\,\sum_{\nu=s-\mu}^{N-\mu}\,
       w_{\mu\nu}\,
        \cdot\, \left({\small \begin{array}{c}
                               \nu \\ s-\mu
                              \end{array}}\right)\,
          d^{s-\mu}\:\lambda_i^{\mu+\nu-s}\;
    +\; \sum_{\nu=0}^{N-s}\,
         w_{s\nu}\,\cdot\,\left(\rule{0em}{1.2em}1+(-1)^{s+1}\right)
          \lambda_i^{\nu}\;
       =\;0  
  $$
  \hspace{4em} for $i=0,\,\cdots,\,n$ and $s=0,\,\cdots,\,N$.
  } 
\end{quote}

\bigskip

\noindent $\bullet$
{\it System from the special values}$\,$:

\begin{quote}
 \hspace{-1.9em}(3)\hspace{1ex} {\it 
  $Q_d(\lambda_i,\frac{\lambda_i-\lambda_j}{d})$ for
  $i=0,\,\cdots,\,n$}$\,$:
  \medskip

  Fact 2.6 gives rise to the third linear system in $w_{\mu\nu}\,$: 

  \raisebox{-.6em}{$\left\{ \rule{0em}{3em}\right.$} \hspace{1em}
  \parbox{12cm}{
   $$
    \sum_{\mu=0}^N\,\sum_{\nu=0}^{N-\mu}\,
       w_{\mu\nu}\,\cdot\,
        \left(\, \frac{\lambda_i-\lambda_j}{d}\,\right)^{\mu}\,
        \lambda_i^{\nu}\;
    =\; \prod_a\prod_{m=0}^{l_ad}\,
       (\,l_a\lambda_i-m\frac{\lambda_i-\lambda_j}{d}\,) \hspace{6em}
   $$
   \hspace{4em} for $i,j=0,\,\cdots,\,n$, $i\ne j$.
  } 
\end{quote}

\bigskip

\noindent $\bullet$
{\it System from the $\alpha$-degree bound}$\,$:

\begin{quote}
 \hspace{-1.9em}(4)\hspace{1ex} {\it
  $\degree_{\alpha}I^{\ast}_d(Q_d)\le (n+1)d-2$ }$\,$:
  \medskip

  Since
  $H^{\ast}_G(\CP^n)\,
    =\, \left. {\Bbb C}[\lambda, \alpha][\kappa]
        \mbox{\raisebox{-2ex}{\rule{0ex}{1ex}}} \right/\hspace{-.4ex}
      \mbox{\raisebox{-.4ex}{$\prod_{i=0}^n(\kappa-\lambda_i)$}}$,
  $I^{\ast}_d(Q_d)$ is obtained by $Q_d$ modulo the relation
  $\prod_{i=0}^n(\kappa-\lambda_i)=0$. This is achieved by iterations
  of the set of replacements
  $$
   \{\, \kappa^{n+1+i}\rightarrow
       \kappa^i\,
        (-\kappa^{n+1}+\mbox{$\prod_{i=0}^n\,(\kappa-\lambda_i)$})\;
        |\, i=0,\, \cdots,\, N-n-1\,\}
  $$
  until the $\kappa$-degree of the resulting $Q_d$ is less than $n+1$.
  This can be easily done by computer. In this way, one obtains
  $I^{\ast}_d(Q_d)$. Let
  $$
   I^{\ast}_d(Q_d)\;=\;\sum_{i,j}\,w^{\prime}_{ij}\alpha^i\kappa^j\,.
  $$
  Then $w^{\prime}_{ij}$ is a linear combination of $w_{\mu\nu}$ with
  coefficients in ${\Bbb C}[\lambda]$. The $\alpha$-degree bound then
  gives us the fourth system of linear equations in $w_{\mu\nu}\,$:

  \raisebox{-.2em}{$\left\{ \rule{0em}{1em}\right.$} \hspace{1em}
  \parbox{12cm}{
   $$
    w^{\prime}_{ij}\;=\;0,
    \hspace{2em}\mbox{for}\hspace{1ex} i\,\ge\,(n+1)d-1\,.
    \hspace{16em}
   $$
  } 
\end{quote}

\bigskip

\noindent
{\it Remark 3.1.}
Note that the whole content of gluing identities in the definition
of Euler data is already absorbed in the first and the second linear
systems (1) and (2) above. Since all the identities that appear in
the Reciprocity Lemma are obtained by substituting into the gluing
identities some special $\alpha$ values, they will be automatically
satisfied once System (1) and System (2) above are satisfied. Thus,
they do not provide us with extra equations.
Furthermore, Fact 2.6 implies that the above system has a unique
solution.

\bigskip

\section{Modifications for non-critical bundles over
         ${\LargeBbb C}{\rm P}^n$.}

So far, our discussion has been focusing on critical bundles.
For non-critical bundles $V$ over $\CP^n$, $\rank(U_d)$ and
$\mbox{\it dim}\,\overline{\cal M}_{0,0}(\CP^n, d)$ are different
for some $d$, by definition. If one takes $\{Q_d\}_d$ associated 
to the top Chern class of $V$, then one will simply get $0$ for
$K_d$. Let $r$ be the rank of $V$. Then, to obtain more interesting
invariants, one may consider taking $\{Q_d\}_d$ to be the Euler data
associated to the Chern polynomial
$$
 c(x)\; =\; x^r + c_1(V)\, x^{r-1}\,+\, \cdots,\, +\, c_{top}(V)
$$
of $V$. In this case, some details in Sec.\ 2 and Sec.\ 3 will
have to be modified accordingly; however, the conceptual flow
remains the same.

The various items that are involved in the actual computation of
$Q_d$ and $K_d$ and their modifications following [L-L-Y] are
listed below:
\begin{itemize}
 \item
 {\bf degree bound of $Q_d$}$\,$:
 Since the top Chern class gives the highest degree terms, both
 the total degree bound and the $\alpha$-degree bound remain valid.
 
 \item
 {\bf the special values
  $\left. i^{\ast}_{p_{i,0}}(Q_d)
    \right|_{\alpha=\frac{\lambda_i-\lambda_j}{d}}$}$\,$:
  These special values are now given by 
  $$
   \left. i^{\ast}_{p_{i,0}}(Q_d)
    \right|_{\alpha=\frac{\lambda_i-\lambda_j}{d}}\;
   =\; \prod_a\prod_{m=0}^{l_ad}\,
             (x+l_a\lambda_i-m\,\frac{\lambda_i-\lambda_j}{d})\:
       \prod_b\prod_{m=1}^{k_bd-1}\,
             (x-k_b\lambda_i+m\,\frac{\lambda_i-\lambda_j}{d})\,.
  $$
  This leads to a corresponding change in the computation of $Q_1$
  and right-hand-side of the third linear systems in Sec.\ 3. The
  first, second, and the fourth linear systems remains valid.

  \item
  {\bf from $Q_d$ to $K_d$}$\,$:
  Let $s=\rank(U_d)-\mbox{\it dim}\,\overline{\cal M}(\CP^n, d)$.
  Theorem 7.2 in [L-L-Y, II] gives 
  $$
   \frac{1}{s!}\,\left.\frac{d^s}{dx^s}\right|_{x=0}
    \int_{{\scriptsizeBbb C}{\rm P}^n}
     e^{-Ht/{\alpha}}\,
       \frac{ \lim_{\lambda\rightarrow 0} I^{\ast}_d(Q_d)}{
               \prod_{m=1}^d(H-m\alpha)^{n+1}}\;
   =\; \frac{1}{\alpha^3\,x^s}\,(2-dt)\,K_d\,.
  $$
  Note that, for the validity of this formula alone, it is not
  required that $s_d$ be independent of $d$ [Li]. Thus, once
  $Q_d$ is obtained by solving the system of linear equation,
  $K_d$ follows.
\end{itemize}

Based on Sec.\ 3 and the discussion here, a code
{\tt \verb+EulerData_MP.m+} is written. The detail is in Sec.\ 6.

\bigskip

\section{Examples.}

Using the Maple code "{\tt \verb+EulerData_MP.m+}" in Sec.\ 6 and
taking $\Omega$ to be the Chern polynomial, we compute the first few
$K_d$ for some non-critical bundles, as listed in {\sc Table 5-1}
(cf.\ Cases 6, 14-18 in SEC.\ 3 of the code).

\bigskip

\noindent\hspace{-1.4em}
{\scriptsize
\begin{tabular}{|l||crrrc|} \hline
 \rule{0ex}{3ex} bundle &  & $d=1$\hspace{1em}
   & $d=2$\hspace{1em} & $d=3$\hspace{1em} &  \\[.6ex] \hline\hline
 $T_{\ast}{\scriptsizeBbb C}{\rm P}^2$ \rule{0ex}{3.2ex}
     & & $10\, x^2$         &  ---  &  ---  &   \\
 $O_{{\tinyBbb C}{\rm P}^4}(6)$ \rule{0ex}{3.2ex}
     & & $50400\, x$  &  $(752729895/4)\, x^2$
                &  $(433244745198080/243)\, x^3$  & \\
 $O_{{\tinyBbb C}{\rm P}^4}(7)$ \rule{0ex}{3.2ex}
     & & $451570\, x^2$
                & $(403985396325/32)\, x^4$
                & $(15755269694706695755/17496)\, x^6$ & \\
 $O_{{\tinyBbb C}{\rm P}^4}(8)$  \rule{0ex}{3.2ex}
     & & $2773820\, x^3$
                & $(3178734062035/8)\, x^6$
                & $(46028387589557254161275/314928)\, x^9$&  \\
 $O_{{\tinyBbb C}{\rm P}^4}(9)$  \rule{0ex}{3.2ex}
     & & $13198850\, x^4$
                & $(243281907041715/32)\, x^8$
                & $(197802281929974511821535/17496)\, x^{12}$ &  \\
 $O_{{\tinyBbb C}{\rm P}^4}(10)$  \rule{0ex}{3.2ex}
     & & $52040450\, x^5$
                & $(25908993204089625/256)\, x^{10}$
                & $(71418501571607082433686025/139968)\, x^{15}$
                                               &  \\[.6ex]  \hline
\end{tabular}
} 

\medskip

\begin{center}
 \parbox{11cm}{{\sc Table 5-1}.
 $K_d$ for some non-critical bundles.
 }
\end{center}

\bigskip

Note that, in {\sc Table 5-1}, if one formally converts $K_d$ to $n_d$
using the usual multiple cover formula
$K_d = \sum_{k|d}\, n_{\frac{d}{k}}\frac{1}{k^3}$, then the $n_d$ thus
obtained will no longer be integers for $d\ge 2$. Compared with the
other twelve examples computed/tested using the code, this indicates
that for those non-critical bundles over $\CP^n$, whose
$\mbox{\it rank}\,(U_d)
           -\mbox{\it dim}\,(\overline{\cal M}_{0,0}(\CP^n, d))$
depends on $d$ in a non-trivial way, the above multiple cover formula
will have to be modified. Exactly how is an issue for further
investigation.

\bigskip

\section{A package in Maple V for the computation of $Q_d$ and $K_d$.}

\begin{flushleft}
{\bf General remarks.}
\end{flushleft}
A few remarks are given below concerning the code, its current scope,
and its usage.

\begin{itemize}
 \item
 Two bad things first. The first one is that the code, as currently
 written, is limited only to the case that $\Omega$ is an Euler class
 or a Chern polynomial. For other multiplicative characteristic
 classes, one will not be able to use the code to make sensible
 computation without entering the core part,
 SEC.\ 1.3 (Q1K1) and SEC.\ 2 (EulerData), to do some modifications.
 This shortage will hopefully be gradually removed along with the
 development of the theory, using the relation of a given
 multiplicative class with Chern roots.
 The second one is that the actual computation of Euler data and
 $K_d$ could be a very demanding task both for Maple V and the
 machine 
 (cf.\ some words in SEC.\ 3 of the code about
                        Test 6: $T_{\ast}\CP^2$ even for $d=2$).
 When it exceeds the capacity the Maple V, one will likely get
 an error message. Experience tells us that occasionally these
 messages may be misleading.
 
 \item
 Now let us mention something more positive.
 As indicated in SEC.\ 3 of the code, the code has been tested
 correct for all known cases within the capacity of Maple V and
 the machine used. It is also made user-friendly:
 {\it to run for a case of study, one only has to follow the
 examples in SEC.\ 3 and modify the various arguments/parameters to
 be fed into the function `$\,${\tt EulerData}', as instructed there.}
 By no means does one need to do anything else; {\it nor} do we
 assume any knowledge of Maple at all.

 \item
 The specialization used, the time consumed, and the RAM memory used
 for the examples tested, particularly those that take long hours,
 are recorded in SEC.\ 3 of the code for reference. Only the last
 run of each case is recorded in these notes.
\end{itemize}

\bigskip

\begin{flushleft}
{\bf Instruction of running the code under Window 98.}
\end{flushleft}
For non-maple-user, let us give here some instruction of running
the code under Window 98: (assuming there is already Maple V in the
folder ``{\tt Programm Files}")

\begin{itemize}
 \item
 Read first the instructions both at the beginning of the code and
 at the start of SEC.\ 3 of the code .

 \item
 Save the code (say by the file name `\verb+EulerData_MP.m+')
 in the subfolder `{\tt Bin.wnt}' in the folder
 `{\tt Maple V Release 5}' that is automatically set up in the
 folder `{\tt Programm Files}' when installing Maple V. Double
 click the maple icon to open a Maple V worksheet. Inside the
 worksheet, type in the following command line after the prompt
 \begin{quote}
 {\tt \verb+ restart; read "EulerData_MP.m"; +}
 \end{quote}
 then hit \fbox{$\,${\small Enter}$\,$} on the keyboard.
 The output will be displayed directly on the worksheet.
\end{itemize}

For all other operating systems, please consult the system manager. 

\newpage

\begin{flushleft}
{\bf The code `{\tt \verb+EulerData_MP.m+}'.}
\end{flushleft}

{\scriptsize\tt


} 

\newpage

\enddocument